\documentclass{article}

\usepackage{amsthm,amsmath,amssymb}
\usepackage{graphicx}
\usepackage{color}

\newtheorem{theorem}{Theorem}
\newtheorem{lemma}{Lemma}

\newtheorem{definition}{Definition}
\newtheorem{proposition}{Proposition}
\newtheorem{remark}{Remark}
\newtheorem{example}{Example}

\newcommand{\R}{\ensuremath{\mathbb R}}

\newcommand{\ls}[1]
    {\dimen0=\fontdimen6\the\font\lineskip=#1\dimen0
     \advance\lineskip.5\fontdimen5\the\font
     \advance\lineskip-\dimen0
     \lineskiplimit=0.9\lineskip
     \baselineskip=\lineskip
     \advance\baselineskip\dimen0
     \normallineskip\lineskip\normallineskiplimit\lineskiplimit
     \normalbaselineskip\baselineskip
     \ignorespaces}

\begin{document}

\bibliographystyle{abbrv}

\title{The fundamental theorem of affine geometry in $(L^0)^n$ }
\author{Mingzhi Wu$^1$ \quad Long Long$^{2}$\\
1. School of Mathematics and Physics, China University of Geosciences,\\ Wuhan {\rm 430074}, China \\
Email: wumz@cug.edu.cn\\
2. School of Mathematics and Statistics,
Central South University,\\
Changsha {\rm 410083}, China \\
Email: longlong@csu.edu.cn\\}

\date{}
 \maketitle

\thispagestyle{plain}
\setcounter{page}{1}

\begin{abstract}

Let $L^0$ be the algebra of equivalence classes of real valued random variables on a given probability space, and $(L^0)^n$ the $n$-ary Cartesian power of $L^0$ for each integer $n\geq 2$. We consider $(L^0)^n$ as a free module over $L^0$ and study affine geometry in $(L^0)^n$. One of our main results states that: an injective mapping $T: (L^0)^n\to (L^0)^n$ which is local and maps each $L^0$-line onto an $L^0$-line must be an $L^0$-affine linear mapping. The other main result states that: a bijective mapping $T: (L^0)^n\to (L^0)^n$ which is local and maps each $L^0$-line segment onto an $L^0$-line segment must be an $L^0$-affine linear mapping. These results extend the fundamental theorem of affine geometry from $\mathbb R^n$ to $(L^0)^n$.
\end{abstract}



\ls{1.5}
\section{Introduction}
The fundamental theorem of affine geometry is a classical and useful result. It states that for an integer $n\geq 2$, if a bijective mapping $F: {\mathbb R}^n\to {\mathbb R}^n$ maps any line to a line, then it must be affine linear. A short proof can be found in Remark 6 of Artstein-Avidan and Milman \cite{AM}.

The fundamental theorem of affine geometry has been generalized and strengthened in numerous ways. Please see Section 5 of Artstein-Avidan and Slomka \cite{AS} for an account of the various forms and generalizations
of the fundamental theorems of affine geometry for ${\R}^n$, together with references and other
historical remarks, and see Kvirikashvili and Lashkhi \cite{KL} and the references therein for generalizations of the fundamental theorems of affine geometry for free modules over some kinds of rings and other more general underlying structures.

Let $L^0$ be the algebra of equivalence classes of real valued random variables on a given probability space, and for any positive integer $n$,  $(L^0)^n=\{(\xi_1,\dots,\xi_n): \xi_i\in L^0,i=1,\dots,n\}$. Since $L^0$ and $(L^0)^n$ (endowed with the usual topology of convergence in probability) usually fail to be local convex spaces, most of mathematicians are not interested in $L^0$ and $(L^0)^n$ for quite a
long time. However, motivated by financial applications and stochastic optimizations, the study of $L^0$ and $(L^0)^n$ became active in the literature recently. For example, Kardaras \cite{Kar} studied the uniform integrability and the local convexity in $L^{0}$; \v{Z}itkovi\'{c}\cite{Zit}, Kardaras and \v{Z}itkovi\'{c}\cite{KZ} considered the forward convex convergence in $L^0$'s nonnegative orthant $L^0_+$; Drapeau, et. al \cite{DKK} established the Brouwer fixed point theorem in $(L^0)^n$; Wu \cite{Wu} established the Farkas' lemma and Minkowski-Weyl type results in $(L^0)^n$ and Cheridito, et. al \cite{CKV} generalized some classical results from linear algebra, real analysis and
convex analysis to $(L^0)^n$.

In this paper, for each $n\geq 2$, we consider $(L^0)^n$ as a free $L^0$-module of rank $n$ and study affine geometry in $(L^0)^n$. We extend the fundamental theorem of affine geometry from $\mathbb R^n$ to $(L^0)^n$. One of our main results states that: an injective mapping $T: (L^0)^n\to (L^0)^n$ which is local and maps each $L^0$-line onto an $L^0$-line must be an $L^0$-affine linear mapping. The other main result states that: a bijective mapping $T: (L^0)^n\to (L^0)^n$ which is local and maps each $L^0$-line segment onto an $L^0$-line segment must be an $L^0$-affine linear mapping. Besides, we also give an example to show that a bijective mapping $T: (L^0)^n\to (L^0)^n$ which maps $L^0$-lines onto $L^0$-lines is not necessarily $L^0$-affine linear.

\section{Basic notations and definitions}

Throughout this paper, $(\Omega,{\mathcal F},P)$ will denote a given probability space. Let $L^0$ be the set of all equivalence classes of real valued
random variables on $(\Omega,{\mathcal F},P)$. Under the usual addition and multiplication operations, $L^0$ is an algebra. For each $A\in {\mathcal F}$, the equivalence class of $A$ refers to $\tilde A=\{B\in {\mathcal F}: P(A\triangle B)=0\}$, ${\tilde I}_A$ denotes the equivalence class of the characteristic function $I_A$, and $I_{\tilde A}$ also stands for ${\tilde I}_A$. Given $\xi\in L^0$, let $\xi^0(\cdot)$ be an arbitrarily chosen representative of $\xi$. We write $[\xi\neq 0]$ for the equivalence class of the measurable set $\{\omega\in\Omega: \xi^0(\omega)\neq 0\}$. The statement ``$\xi\neq 0$ on $\Omega$'' means that $\xi^0(\omega)\neq 0$, $P$-a.s., in other words, $\xi$ is an invertible element of the algebra $L^0$. Some other notation like $[\xi=0]$ or statement like ``$\xi>0$ on $\Omega$'' and so on are understood in a similar way.

In this paper, an $L^0$-module refers to a left module over the algebra $L^0$, and $\theta$ always denotes its null element.

Given any positive integer $n$, denote $(L^0)^n=\{(\xi_1,\dots,\xi_n): \xi_i\in L^0,i=1,\dots,n\}$, then $(L^0)^n$ is a free $L^0$-module of rank $n$ generated by $e_i,i=1,\dots,n$, where $e_i$ is the $i$-th unit vector in ${\mathbb R}^n\subset (L^0)^n$.  Given $x=(\xi_1,\dots,\xi_n), y=(\eta_1,\dots,\eta_n)\in (L^0)^n$, the $L^0$-inner product of $x$ and $y$ is defined to be $\langle x,y\rangle=\sum^n_{i=1}\xi_i\eta_i$, and the $L^0$-normed of $x$ is given by $|x|=\sqrt{\langle x,x\rangle}$. We can see that $|x|=0$ iff $x=\theta$ and $|\xi x|=|\xi||x|$ for every $x\in (L^0)^n$ and $\xi\in L^0$.

Given $x\in (L^0)^n$, if $|x|\neq 0$ on $\Omega$, then $x$ is said to have full support. Obviously, each $e_i$ has full support. If $x\in (L^0)^n$ has full support, then for any $\xi\in L^0$, using the equality $|\xi x|=|\xi||x|$, we see that $\xi x=\theta$ iff $\xi=0$; as a result, if $\xi$ and $\eta$ are elements of $L^0$ such that $\xi x=\eta x$, then $\xi=\eta$. On the contrary, if $x\in (L^0)^n$ does not have full support, then $A=[|x|=0]$ has positive probability and $I_Ax=\theta$, immediately $x=I_{A^c}x$, where $A^c=[|x|\neq 0]$ so that $I_{A^c}=1-I_A$, and it is probably that $I_{A^c}\neq 1$.

 A group of elements $x_1,\dots,x_k$ in $(L^0)^n$ are said to be $L^0$-independent, if for $\xi_1,\dots,\xi_k \in L^0$, the equality $\xi_1 x_1+\cdots+\xi_k x_k =\theta$ implies that $\xi_1=\cdots=\xi_k=0$. Since $I_{[|x_1|=0]}x_1+\cdots+I_{[|x_k|=0]}x_k=\theta+\cdots+\theta=\theta$, we see that all elements $x_i$'s have full support whenever $x_1,\dots,x_k$ are $L^0$-independent. Clearly, $e_1,\dots,e_n$ are $L^0$-independent.

Fix an integer $n\geq 2$. For any nonzero vector $x\in \mathbb R^n$, we can find a vector $y\in \mathbb R^n$ linearly independent of $x$, while for a nonzero element $x\in (L^0)^n$, we may not find an element $y\in (L^0)^n$ such that $y$ and $x$ are $L^0$-independent, in fact, according to the aforementioned fact, for a nonzero $x$ which does not have full support, there does not exist $y$ which is $L^0$-independent of $x$ at all! When $x\in (L^0)^n$ has full support, the existence of $y\in (L^0)^n$ which is $L^0$-independent of $x$ is also not obvious. These observations force us to use the base $\{e_1,\dots,e_n\}$ in the proof of our main result Theorem \ref{main}.

\begin{definition}\label{def}
Let $E_1$ and $E_2$ be two $L^0$-modules and $T: E_1 \to E_2$ a mapping.\\
(1). $T$ is said to be $L^0$-linear, if $T(x+y)=T(x)+T(y), \forall x,y\in E_1$, and $T(\xi x)=\xi T(x), \forall x\in E_1, \xi\in L^0$;\\
(2). $T$ is said to be $L^0$-affine linear, if $T(\lambda x+(1-\lambda)y)=\lambda T(x)+(1-\lambda)T(y), \forall x,y\in E_1, \lambda\in L^0$, equivalently, $T(\cdot)-T(\theta)$ is $L^0$-linear; \\
(3). $T$ is said to be local (or have the local property) if ${\tilde I}_AT({{\tilde I}_A x})={\tilde I}_AT(x), \forall x\in E_1, A\in {\mathcal F}$;\\
(4). $T$ is said to be stable if for any $x,y\in E_1$ and $A\in {\mathcal F}$, $T({\tilde I}_A x+{\tilde I}_{A^c}y)={\tilde I}_AT(x)+{\tilde I}_{A^c}T(y)$.
\end{definition}

\begin{proposition}\label{aff-lc}
 Let $E_1$ and $E_2$ be two $L^0$-modules and $T: E_1 \to E_2$ a mapping. Then the following statements are true: \\
 (1). If $T$ is $L^0$-affine linear, then $T$ must have the local property.\\
 (2). $T$ has the local property if and only if $T$ is stable.\\
 (3). If $T$ has the local property and $T(\theta)=\theta$, then $T({\tilde I}_A x)={\tilde I}_A T(x)$ for any $x\in E_1$ and $A\in {\mathcal F}$.\\
 (4). If $T$ is bijective and has the local property, then $T^{-1}$ also has local property.
\end{proposition}

\begin{proof}
 (1). Define the $L^0$-linear mapping $S: E_1\to E_2$ by $S(x)=T(x)-T(\theta), \forall x\in E_1$. For any $x\in E$ and $A\in {\mathcal F}$, ${\tilde I}_AT({\tilde I}_Ax)={\tilde I}_A[S({\tilde I}_Ax)+T(\theta)]={\tilde I}_A[{\tilde I}_AS(x)+T(\theta)]={\tilde I}_A[S(x)+T(\theta)]={\tilde I}_AT(x)$, namely, $T$ has the local property.

(2). If $T$ has the local property, then for any $x,y\in E_1$ and $A\in {\mathcal F}$
\begin{eqnarray*}
T({\tilde I}_A x+{\tilde I}_{A^c}y)&=& {\tilde I}_AT({\tilde I}_A x+{\tilde I}_{A^c}y)+ {\tilde I}_{A^c}T({\tilde I}_A x+{\tilde I}_{A^c}y) \\
                                   &=& {\tilde I}_AT({\tilde I}_A({\tilde I}_A x+{\tilde I}_{A^c}y))+{\tilde I}_{A^c}T({\tilde I}_{A^c}({\tilde I}_A x+{\tilde I}_{A^c}y))\\
                       & =& {\tilde I}_AT({\tilde I}_A x)+ {\tilde I}_{A^c}T({\tilde I}_{A^c}y) \\
                       & =& {\tilde I}_AT(x)+{\tilde I}_{A^c}T(y).
\end{eqnarray*}
Thus, $T$ is stable.

Conversely, if $T$ is stable, then for any $x\in E_1$ and $A\in {\mathcal F}$, $T({\tilde I}_A x)=T({\tilde I}_A x+{\tilde I}_{A^c}\theta)={\tilde I}_AT(x)+ {\tilde I}_{A^c}T(\theta)$, immediately we obtain ${\tilde I}_AT({\tilde I}_A x)={\tilde I}_AT(x)$, which means that $T$ has the local property.

(3). From (2), $T({\tilde I}_A x)=T({\tilde I}_A x+{\tilde I}_{A^c}\theta)={\tilde I}_AT(x)+ {\tilde I}_{A^c}T(\theta)={\tilde I}_AT(x)$.

(4). Fix $A\in {\mathcal F}$ and $y\in E_2$. Using (2), $T[{\tilde I}_A T^{-1}(y)]={\tilde I}_A T[T^{-1}(y)]+{\tilde I}_{A^c}T(\theta)={\tilde I}_Ay+{\tilde I}_{A^c}T(\theta)$, and $T[{\tilde I}_AT^{-1}({\tilde I}_Ay)]={\tilde I}_AT[T^{-1}({\tilde I}_Ay)]+{\tilde I}_{A^c}T(\theta)={\tilde I}_Ay+{\tilde I}_{A^c}T(\theta)$. Then the assumption that $T$ is injective yields that ${\tilde I}_A T^{-1}(y)={\tilde I}_AT^{-1}({\tilde I}_Ay)$, exactly meaning $T^{-1}$ also has the local property.
\end{proof}

\section{Main results}

For any two distinct points $x,y$ in $(L^0)^n$, denote $l(x,y)=\{\lambda x+(1-\lambda)y: \lambda\in L^0\}$, called the $L^0$-line determined by $x$ and $y$. In $\mathbb R^n$, any two distinct points in a given straight line determine the same straight line, while in $(L^0)^n$, if $u,v$ are two distinct points in the $L^0$-line $l(x,y)$, the $L^0$-line $l(u,v)$ may be not the same as $l(x,y)$. For instance, let $x\in (L^0)^n$ be a nonzero element, then for any $A\in {\mathcal F}$, ${\tilde I}_Ax$ lies in the $L^0$-line $l(\theta,x)=\{\lambda x: \lambda\in L^0\}$, if ${\tilde I}_Ax$ is nonzero, then the $L^0$-line $l(\theta,{\tilde I}_Ax)=\{{\tilde I}_A\lambda x: \lambda\in L^0\}$ is probably not the same as $l(\theta,x)$. Thus we should be careful when we handle problems involving $L^0$-lines.

It is easy to verify that any injective $L^0$-affine linear mapping $T: (L^0)^n\to (L^0)^n$ maps each $L^0$-line onto an $L^0$-line. Theorem \ref{main} below states that the converse is also true.

\begin{theorem}\label{main}
 Fix an integer $n\geq 2$. Let $T: (L^0)^n\to (L^0)^n $ be an injective mapping which is local and maps each $L^0$-line onto an $L^0$-line, that is to say, for any two distinct points $x,y\in (L^0)^n$, the image of the $L^0$-line $l(x,y)$ under the mapping $T$ is $l(u,v)$, where $u=T(x),v=T(y)$, then $T$ must be an $L^0$-affine linear mapping.
\end{theorem}

\begin{proof}

 Define $S: (L^0)^n\to (L^0)^n$ by $S(x)=T(x)-T(\theta), \forall x\in (L^0)^n$. Note that $S(\theta)=\theta$ and $S$ is also injective. With the assumptions on $T$, it is easy to check that $S$ is local and maps each $L^0$-line onto an $L^0$-line. It remains to show that $S$ is $L^0$-linear. The proof is composed of 5 steps as below. We point out in advance that (2) and (3) of Proposition \ref{aff-lc} are used frequently.

{\em Step 1.} For any $L^0$-independent $x,y\in (L^0)^n$, we have $S(x)$ and $S(y)$ are $L^0$-independent, and $S(x+y)=S(x)+S(y)$.

For any $z\in (L^0)^n$ which has full support, let $A=[|S(z)|=0]$, then by (3) of Proposition \ref{aff-lc}, $S(I_A z)=I_AS(z)=\theta$. Since $S$ is injective, we obtain that $I_A z=\theta$. Thus $I_A=0$, implying $S(z)$ has full support.

Suppose $\xi,\eta\in L^0$ satisfy the equality $\xi S(x)+\eta S(y)=\theta$. Since $S$ is injective and maps the $L^0$-line $l(\theta,x)$ onto the  $L^0$-line $l(\theta, S(x))$, there exists $\alpha\in L^0$ such that $\xi S(x)=S(\alpha x)$. Similarly, there exists $\beta\in L^0$ such that $-\eta S(y)=S(\beta y)$. By the injectivity of $S$ we get $\alpha x=\beta y$, then $\alpha=\beta=0$ follows from the assumption that $x,y$ are $L^0$-independent. As a result, $\xi S(x)=-\eta S(y)=\theta$, then using the fact that both $S(x)$ and $S(y)$ have full support, we conclude that $\xi=\eta=0$, which means that $S(x)$ and $S(y)$ are $L^0$-independent.

We then show: there exist $a,b\in L^0$ such that $S(x+y)=aS(x)+bS(y)$. In fact, since $x+y$ lies in the $L^0$-line $l(2x,2y)$ and $S$ maps $L^0$-lines to $L^0$-lines, thus there exists $\mu\in L^0$ such that $S(x+y)=\mu S(2x)+(1-\mu)S(2y)$. Since $2x$ lies in the $L^0$-line $l(\theta,x)$ and $2y$ lines in the $L^0$-line $l(\theta,y)$, there exist $\alpha_1,\beta_1\in L^0$ such that $S(2 x)=\alpha_1S(x), S(2y)=\beta_1S(y)$, then $a=\mu\alpha_1,b=(1-\mu)\beta_1$ satisfy the equality $S(x+y)=aS(x)+bS(y)$.

 It remains to show $a=1$ and $b=1$. Let $A=[a-1\neq 0]$, if by contrary that $a\neq 1$, then $A$ has positive probability and $I_A\neq 0$. According to the notation, there exists $c_1\in L^0$ such that $I_A[1+c_1(a-1)]=0$. Since the $L^0$-line $l(x,x+y)=\{x+cy: c\in L^0\}$ is mapped by $S$ onto the $L^0$-line $l(S(x),S(x+y))$, there exists $c_0\in L^0$ such that $S(x+c_0y)=(1-c_1)S(x)+c_1S(x+y)=[1+c_1(a-1)]S(x)+c_1bS(y)$. Using Proposition \ref{aff-lc} we obtain $S(I_A(x+c_0y))=I_AS(x+c_0y)=I_Ac_1b S(y)$. Note that there exists some $\xi\in L^0$ such that ${\tilde I}_Ac_1b S(y)=S(\xi y)$, then by the injectivity of $S$, we get $I_A(x+c_0y)=\xi y$, contradicting to the assumption that $x,y$ are $L^0$-independent. Therefore, $a=1$. Similarly, $b=1$.

{\em Step 2.} For any $L^0$-independent $x,y\in (L^0)^n$, we have $S(\xi x+\eta y)=S(\xi x)+S(\eta y), \forall\xi,\eta\in L^0$.

First suppose $\xi,\eta$ are characteristic functions, that is $\xi=\tilde I_A$, $\eta=\tilde I_B$, for some $A,B\in {\mathcal F}$. Since
$\tilde I_A x+\tilde I_By=\tilde I_{A\cap B}(x+y)+\tilde I_{A\setminus B}x+\tilde I_{B\setminus A}y$, by the local property, we have $S(\tilde I_A x+\tilde I_By)=\tilde I_{A\cap B}S(x+y)+\tilde I_{A\setminus B}S(x)+\tilde I_{B\setminus A}S(y)=\tilde I_{A\cap B}[S(x)+S(y)]+\tilde I_{A\setminus B}S(x)+\tilde I_{B\setminus A}S(y)=\tilde I_A S(x)+\tilde I_BS(y)=S(\tilde I_A x)+S(\tilde I_By)$.

Generally, for any $\xi,\eta\in L^0$, let $A=[\xi\neq 0], B=[\eta\neq 0]$, and take $x_1=\xi x+I_{A^c}x, y_1=\eta y+I_{B^c}y$, then $\xi x=I_Ax_1, \eta y=I_B y_1$, and $x_1,y_1$ are $L^0$-independent. In fact, if $\alpha,\beta\in L^0$ satisfy the equality $\alpha x_1+\beta y_1=\theta$, then since $x,y$ are $L^0$-independent, we must have $\alpha(\xi+I_{A^c})=0$ and $\beta(\eta+I_{B^c})=0$, due to the fact $\xi+I_{A^c}\neq 0$ on $\Omega$ and $\eta+I_{B^c}\neq 0$ on $\Omega$, we thus obtain $\alpha=\beta=0$. Now we have shown that $x_1,y_1$ are $L^0$-independent, then $S(\xi x+\eta y)=S(I_A x_1+I_By_1)=S(\tilde I_A x_1)+S(\tilde I_By_1)=S(\xi x)+S(\eta y)$.

{\em Step 3.} For each $i\in \{1,2\dots,n\}$ and any $\xi,\eta\in L^0$, we have $S(\xi e_i+\eta e_i)=S(\xi e_i)+S(\eta e_i)$.

By symmetry, it suffices to prove the case when $i=1$.

Since $e_1-e_2$ and $e_2$ is obvious $L^0$-independent, we get from Step 2 that $S(e_1)=S(e_1-e_2+ e_2)=S(e_1-e_2)+S(e_2)=S(e_1)+S(-e_2)+S(e_2)$, therefore $S(e_2)+S(-e_2)=\theta$.

Now fix $\xi,\eta\in L^0$, let $A=[\xi+\eta\neq 0]$. Then $x_1=\xi e_1+I_{A^c}e_1+e_2$ and $y_1=\eta e_1-e_2$ are $L^0$-independent. Indeed, if $\alpha,\beta\in L^0$ satisfy $\alpha x_1+\beta y_1=(\alpha \xi+\alpha I_{A^c}+\beta\eta)e_1+(\alpha-\beta)e_2=\theta$, then $\alpha \xi+\alpha I_{A^c}+\beta\eta=0$ and $\alpha-\beta=0$, equivalently, $\alpha=\beta$ and $\alpha(\xi+\eta+I_{A^c})=0$, thus $\alpha=\beta=0$ follows from the fact that $\xi+\eta+I_{A^c}\neq 0$ on $\Omega$. From Step 2, noting that $e_1$ and $I_{A^c}e_1+e_2$ are $L^0$-independent, we get $S(x_1)=S(\xi e_1)+S(I_{A^c}e_1+e_2)=S(\xi e_1)+S(I_{A^c}e_1)+S(e_2)$ and $S(y_1)=S(\eta e_1)+S(-e_2)$. On the other hand, due to the fact that $x_1$ and $y_1$ are $L^0$-independent, it follows from Step 2 that $S(\xi x+I_{A^c} x+\eta x)=S(x_1+y_1)=S(x_1)+S(y_1)$. Therefore, using the known fact $S(e_2)+S(-e_2)=\theta$ we get $S(\xi e_1+I_{A^c} e_1+\eta e_1)=S(\xi e_1)+S(I_{A^c}e_1)+S(\eta e_1)$. Using the local property, $I_{A^c}S(\xi e_1+I_{A^c} e_1+\eta e_1)=S[I_{A^c}(\xi e_1+I_{A^c} e_1+\eta e_1)]=S(I_{A^c} e_1)$,
hence $S(\xi e_1+\eta e_1)=S[I_A(\xi e_1+I_{A^c} e_1+\eta e_1]=I_AS(\xi e_1+I_{A^c} e_1+\eta e_1)=S(\xi e_1+I_{A^c} e_1+\eta e_1)-I_{A^c}S(\xi e_1+I_{A^c} e_1+\eta e_1)=S(\xi e_1)+S(\eta e_1)$.

 {\em Step 4.} For any $\xi_1,\xi_2,\dots,\xi_n\in L^0$, we have $S(\xi_1e_1+\cdots+\xi_ne_n)=S(\xi_1e_1)+\cdots+S(\xi_ne_n)$.

Indeed, let $y_1=\xi_2e_2+\cdots+\xi_ne_n$ and set $A=[|y_1|\neq 0]$, further take $y=I_A y_1+I_{A^c}e_2$, then $e_1, y$ are $L^0$-independent, and $y_1=I_Ay$, thus according to Step 2, $S(\xi_1e_1+\cdots+\xi_ne_n)=S(\xi_1e_1+y_1)=S(\xi_1e_1+I_Ay)=S(\xi_1e_1)+S(I_Ay)=S(\xi_1e_1)+S(\xi_2e_2+\cdots+\xi_ne_n)$.
By induction, we obtain $S(\xi_1e_1+\cdots+\xi_ne_n)=S(\xi_1e_1)+\cdots+S(\xi_ne_n)$.

{\em Step 5}. For each $i\in \{1,2\dots,n\}$ and any $\xi\in L^0$, we have $S(\xi e_i)=\xi S(e_i)$.

 Fix an $x\in (L^0)^n$ which has full support. Since $\xi x$ lies in the $L^0$-line $l(\theta,x)$, there exists $\mu \in L^0$ such that $S(\xi x)=\mu S(x)$. By Step 1, $S(x)$ has full support, thus this $\mu$ is unique determined by $\xi$ (and $x$).
Therefore, we can define a mapping $f_x:L^0\to L^0$ by the relation $S(\xi x)=f_x(\xi)S(x), \forall \xi \in L^0$.

Specially, for each $i\in \{1,2\dots,n\}$, we have a mapping $f_i: L^0\to L^0$ such that $S(\xi e_i)=f_i(\xi) S(e_i), \forall \xi \in L^0$.

 We show that all mappings $f_i$ are indeed the same one.

 In fact, by symmetry, it suffices to verify that $f_1=f_2$. For each $\xi\in L^0$, on one hand, since $e_1+e_2$ has full support, we have $S(\xi(e_1+e_2))=f_{e_1+e_2}(\xi)S(e_1+e_2)=f_{e_1+e_2}(\xi)[S(e_1)+S(e_2)]$, where the last equality follows from Step 1. On the other hand, from Step 2, $S(\xi(e_1+e_2))=S(\xi e_1)+S(\xi e_2)=f_1(\xi)S(e_1)+f_2(\xi)S(e_2)$. Thus we obtain $f_{e_1+e_2}(\xi)[S(e_1)+S(e_2)]=f_1(\xi)S(e_1)+f_2(\xi)S(e_2)$. Since we have known from Step 1 that $S(e_1)$ and $S(e_2)$ are $L^0$-independent, thus $f_1(\xi)=f_{e_1+e_2}(\xi)=f_2(\xi)$.

Please note that using a similar argument, for any $\eta\in L^0$ such that $\eta\neq 0$ on $\Omega$, we have $f_{\eta e_1}(\xi)=f_2(\xi)=f_1(\xi), \forall \xi\in L^0$.

We proceed to show that $f_1(\xi)=\xi, \forall \xi\in L^0$.

First, it is obvious that $f_1(0)=0$ and $f_1(1)=1$. Then by the local property of $S$, for any $\xi\in L^0$ and $A\in {\mathcal F}$, $S(\tilde I_A\xi e_1)=\tilde I_AS(\xi e_1)$, we obtain that $f_1(\tilde I_A\xi)=\tilde I_Af_1(\xi)$, which means that $f_1$ is local. By Step (3), for any $\xi,\eta\in L^0$, $S(\xi e_1+\eta e_1)=S(\xi e_1)+S(\eta e_1)$, implying that $f_1(\xi+\eta)=f_1(\xi)+f_1(\eta)$. Finally, for any $\xi,\eta\in L^0$, choose $\eta_1\in L^0$ such that $\eta_1\neq 0$ on $\Omega$ and $\eta=I_A\eta_1$, where $A=[\eta\neq 0]$ (for instance, we can take $\eta_1=I_A\eta+I_{A^c}$), then on one hand $S((\xi\eta)e_1)=f_1(\xi\eta)S(e_1)$, on the other hand, $S((\xi\eta)e_1)=S(\xi I_A\eta_1e_1)=I_AS(\xi \eta_1e_1)=I_Af_{\eta_1e_1}(\xi)S(\eta_1 e_1)=I_Af_1(\xi)f_1(\eta_1)S(e_1)=f_1(\xi)f_1(\eta)S(e_1)$. Noting that $S(e_1)$ has full support, we thus obtain $f_1(\xi\eta)=f_1(\xi)f_1(\eta)$.

To sum up, $f_1$ satisfies all the conditions (1-4) in Lemma \ref{id} below, thus $f_1(\xi)=\xi, \forall \xi\in L^0$.

Combining Step 4 and Step 5, we conclude that $S$ is $L^0$-linear, completing the proof.
\end{proof}

\begin{lemma}\label{id}
Let $\phi: L^0\to L^0$ be a mapping such that:\\
(1). $\phi$ is local;\\
(2). $\phi(\xi+\eta)=\phi(\xi)+\phi(\eta), \forall \xi,\eta\in L^0$;\\
(3). $\phi(\xi\eta)=\phi(\xi)\phi(\eta),\forall \xi,\eta\in L^0$;\\
(4). $\phi(1)=1$.\\
Then $\phi$ is the identity, namely, $\phi(\xi)=\xi, \forall \xi\in L^0$.
 \end{lemma}

\begin{proof}
From (2), $\phi(0)+\phi(1)=\phi(1+0)=\phi(1)$, thus $\phi(0)=0$, then $\phi(\xi-\xi)=\phi(0)=\phi(\xi)+\phi(-\xi)$ yields that $\phi(-\xi)=-\phi(\xi)$ for every $\xi\in L^0$. Since $\phi(1)=1$, it is easy to deduce: for any integer $p$, $\phi(p)=p$ and further for any rational number $r$, $\phi(r)=r$. Now assume $q=\sum^d_{i=1}r_i{\tilde I}_{A_i}$ is a simple function in $L^0$ such that every $r_i$ is a rational number, then by the local property of $\phi$, we obtain: $\phi(q)=\sum^d_{i=1}{\tilde I}_{A_i}\phi(r_i)=\sum^d_{i=1}{\tilde I}_{A_i}r_i=q$.

Let $\xi,\eta$ be two elements in $L^0$ with $\xi \geq \eta$, then from (3) we obtain $\phi(\xi-\eta)=\phi(\sqrt {\xi-\eta}\sqrt {\xi-\eta})=\phi(\sqrt {\xi-\eta})\phi(\sqrt {\xi-\eta})\geq 0$. Since $\phi(\xi-\eta)=\phi(\xi)+\phi(-\eta)=\phi(\xi)-\phi(\eta)$, it follows that $\phi(\xi)\geq \phi(\eta)$, that is to say, $\phi$ is monotonically increasing.

For any $\xi\in L^0$, let $q_-=\sum^d_{i=1}r_i{\tilde I}_{A_i}$ and $q_+=\sum^k_{j=1}t_j{\tilde I}_{B_j}$ be any two simple functions in $L^0$ such that every $r_i$ and $t_j$ are rational numbers and $q_-\leq \xi\leq q_+$, then using the monotonicity of $\phi$, we have $q_-=\phi(q_-)\leq \phi(\xi)\leq \phi(q_+)=q_+$. Taking all such possible $q_-$ and $q_+$, we thus obtain that $\phi(\xi)=\xi$, completing the proof.
\end{proof}

\begin{remark}
The local property appears frequently in the study related to $L^0$ and $(L^0)^n$, for example, it appears in the intermediate value theorem of $L^0$- valued functions (Theorem 1.6 of \cite{GZ}) and the Brouwer fixed point theorem in $(L^0)^n$ (Theorem 2.3 in \cite{DKK}, where the local property appears in a slightly more general form which is equivalent to be stable in Definition \ref{def}).
\end{remark}

In the following, we give an example which shows that a bijective mapping $T: (L^0)^n\to (L^0)^n$ which maps any $L^0$-line onto an $L^0$-line may not have the local property, thus according to Proposition \ref{aff-lc}, this mapping $T$ is not $L^0$-affine linear.

\begin{example}\label{exam}

Let $\theta: (\Omega,{\mathcal F},P)\to (\Omega,{\mathcal F},P)$ be an isomorphism, that is to say, $\theta$ is bijective and both $\theta$ and $\theta^{-1}$ are measure-preserving. Then $\theta$ induces a bijection $\sigma: L^0\to L^0$ through $\xi\mapsto $ the equivalence class of $\xi^0(\theta(\cdot))$, where $\xi^0$ is a representative of $\xi\in L^0$. Further, for each positive integer $n$, $\theta$ induces a bijection $T:(L^0)^n\to (L^0)^n$ through $(\xi_1,\dots,\xi_n)\mapsto (\sigma(\xi_1),\dots,\sigma(\xi_n))$. Then it is straightforward to check that $T$ maps each $L^0$-line onto an $L^0$-line. However, if $\theta$ is not the identity mapping, $T$ is probably not local.

As follows is a more concrete example.

Let $\Omega=[0,1)$, ${\mathcal F}={\mathcal B}([0,1))$, namely the Borel $\sigma$-algebra of $[0,1)$, and $P$ the Lebesgue measure.
Define $\theta: [0,1)\to [0,1)$ by $\theta(\omega)=\omega+\frac{1}{2}$ for $\omega\in [0,\frac{1}{2})$, and $\theta(\omega)=\omega-\frac{1}{2}$ for $\omega\in [\frac{1}{2}, 1)$. We show the induced mapping $T:(L^0)^n\to (L^0)^n$ is not local. Let $A=[0,\frac{1}{2})$, $B=[\frac{1}{2}, 1)$, then $\sigma(\tilde I_A)=\tilde I_B$, therefore for each $x\in (L^0)^n$, we have $T(\tilde I_A x)=\tilde I_B T(x)$, specially $T(\tilde I_A e_1)=\tilde I_B T(e_1)=\tilde I_B e_1$, it follows that $\tilde I_AT(\tilde I_A e_1)=\theta\neq \tilde I_Ae_1=\tilde I_AT(e_1)$. Thus $T$ is not local.

\end{example}

\begin{remark}
 Since $L^0$ is a commutative algebra and $L^0\neq \{0\}$, it follows from Thereom 2.6 in \cite{Cohn} that $L^0$ is an IB-ring (see \cite{KL} for the meaning of this notation), then applying Theorem 1 in \cite{KL} to $(L^0)^n$ gives: for $n\geq 2$, if $T: (L^0)^n\to (L^0)^n$ with $T(\theta)=\theta$ is a collineation preserving parallelism, that is to say, $T$ is a bijection such that the images of collinear points under $T$ are themselves collinear and $T$ preserves parallelism(see \cite{KL} for this notion), then there exists an isomorphism $\sigma: L^0\to L^0$ such that $T$ is a
$\sigma$-semilinear isomorphism, namely $T(x+y)=T(x)+T(y), \forall x,y\in (L^0)^n$ and $T(\xi x)=\sigma(\xi)T(x), \forall x\in (L^0)^n, \xi\in L^0$. Since bijection and preserving parallelism are not premise conditions in our Theorem \ref{main} and our theorem \ref{main} require the local property instead, one can see that our Theorem \ref{main} is not a special case of Theorem 1 in \cite{KL}. We also would like to point out that although $T$ in Example \ref{exam} is not an $L^0$-linear mapping, it is indeed a $\sigma$-semilinear isomorphism.
\end{remark}

For any two distinct $x,y\in (L^0)^n$, denote $[x,y]=\{\mu x+(1-\mu)y: \mu\in L^0, 0\leq \mu\leq 1\}$, called the $L^0$-line segment between $x$ and $y$. In the end of this paper, we discuss self-mappings on $(L^0)^n$ which map $L^0$-line segments to $L^0$-line segments.

\begin{proposition}\label{line-seg}
  Suppose that $T: (L^0)^n\to (L^0)^n $ is a bijection which is local and maps each $L^0$-line segment onto an $L^0$-line segment, that is to say, for any two distinct $x,y\in (L^0)^n$, the image of the $L^0$-line segment $[x,y]$ is the $L^0$-line segment $[Tx,Ty]$, then $T$ maps each $L^0$-line onto an $L^0$-line.
\end{proposition}

\begin{proof}
 Without loss of generality, we can assume $T(\theta)=\theta$, otherwise we make a translation. Let $x,y$ be any two elements in $(L^0)^n$ such that $y\neq \theta$. Since $T$ is a bijection, we can see that $T^{-1}$ also maps each $L^0$-line segment onto an $L^0$-line segment, and $T^{-1}$ is local by Proposition \ref{aff-lc}, thus we only need to show that each point $z$ in the $L^0$-line $l(x,x+y)=\{x+\lambda y: \lambda\in L^0\}$ will be mapped into the $L^0$-line $l(T(x),T(x+y))=\{\lambda T(x)+(1-\lambda)T(x+y): \lambda\in L^0\}$.

We first show that: for each $k\in \mathbb{Z}=\{0,\pm 1,\pm 2,\dots\}$, $z=x+ky$ will be mapped into $l(T(x),T(x+y))$.

First assume that $y$ has full support. (1) The cases $k=0$ and $k=1$ are obvious. (2) Fix a $k\in\{2,3,4,\dots\}$. Since $x+y=(1-\frac{1}{k})x+\frac{1}{k}(x+ky)\in [x,x+ky]$ and $T$ maps an $L^0$-line segment onto an $L^0$-line segment, there exists $\mu\in L^0$ with $0\leq \mu\leq 1$ such that $T(x+y)=(1-\mu) T(x)+\mu T(x+ky)$. Let $A=[\mu=0]$, then $I_A\mu=0$. By the local property and Proposition \ref{aff-lc}, $T(I_A(x+y))=I_AT(x+y)=I_A[(1-\mu) T(x)+\mu T(x+ky)]=I_AT(x)=T(I_Ax)$. Since $T$ is a bijection, we obtain $I_A(x+y)=I_A x$. Then the assumption $y$ has full support implies $I_A=0$, equivalently, $\mu>0$ on $\Omega$. As a result, $T(x+ky)=\frac{1}{\mu}T(x+y)+(1-\frac{1}{\mu})T(x)\in l(T(x), T(x+y))$. (3) Fix a $k\in\{-1,-2,-3,\dots,\}$. Since $x=\frac{1}{1-k}(x+ky)+(1-\frac{1}{1-k})(x+y)\in [x+ky, x+y]$, there exists $\mu\in L^0$ with $0\leq \mu\leq 1$ such that $T(x)=\mu T(x+ky)+(1-\mu)T(x+y)$, by a similar argument we deduce that $\mu>0$ on $\Omega$, then $T(x+ky)=\frac{1}{\mu}T(x)+(1-\frac{1}{\mu})T(x+y)\in l(T(x), T(x+y))$.

 Now for a general nonzero $y$. Take $y_1=I_Ay+I_{A^c}e_1$, where $A=[|y|\neq 0]$, we see that $y_1$ has full support and $y=I_Ay_1$. Fix any $k\in {\mathbb Z}$, we have proved that there exists $\mu \in L^0$ such that $T(x+ky_1)=\mu T(x)+(1-\mu)T(x+y_1)$. Using Proposition \ref{aff-lc},  $T(x+ky)=T[I_A(x+ky)]+T[I_{A^c}(x+ky)]=T[I_A(x+ky_1)]+T(I_{A^c}x)=I_A[\mu T(x)+(1-\mu)T(x+y_1)]+I_{A^c}T(x)=(I_A\mu+I_{A^c})T(x)+I_A(1-\mu)T(x+y)$, implying that $T(x+ky)\in l(T(x),T(x+y))$.

We then show that: for each $\lambda\in L^0$, $z=x+\lambda y$ will be mapped into $l(T(x),T(x+y))$.

For $k=1,2,\dots$, let $A_k=[k-1\leq |\lambda| <k]$, then $x+\lambda I_{A_k}y=(\frac{1}{2}-\frac{\lambda}{2k}I_{A_k})(x-ky)+(\frac{1}{2}+\frac{\lambda}{2k}I_{A_k})(x+ky)$ belongs to $[x-ky, x+ky]$, consequently, $T(x+\lambda I_{A_k}y)\in [T(x-ky), T(x+ky)]\subset l(T(x),T(x+y))$, where the last inclusion follows from Claim 1 that both the two endpoints of the $L^0$-line segment belong to the $L^0$-line. By the notation, for each positive integer $k$, there exists $\mu_k\in L^0$ such that $T(x+\lambda I_{A_k}y)=\mu_kT(x)+(1-\mu_k)T(x+y)$. By the local property of $T$, we have $I_{A_k}T(x+\lambda y)=I_{A_k}T(x+\lambda I_{A_k}y)=I_{A_k}[\mu_kT(x)+(1-\mu_k)T(x+y)]$ for each $k$.
Let $\mu=\sum^{\infty}_{k=1}I_{A_k}\mu_k$, then $T(x+\lambda y)=\sum^{\infty}_{k=1}I_{A_k}T(x+\lambda y)=\sum^{\infty}_{k=1}I_{A_k}[\mu_k T(x)+(1-\mu_k)T(x+y)]=\mu T(x)+(1-\mu)T(x+y)$, which means that $T(x+\lambda y)\in l(Tx,T(x+y))$.
\end{proof}

\begin{remark}
In the proof of Proposition \ref{line-seg}, the local property plays an important role, we wonder whether the assumption that $T$ has the local property can be removed or not? That is,  if $T: (L^0)^n\to (L^0)^n $ is a bijection and maps each $L^0$-line segment onto an $L^0$-line segment, then can we deduce that $T$ maps each $L^0$-line onto an $L^0$-line?
\end{remark}

Combining Theorem \ref{main} and Proposition \ref{line-seg}, we immediately obtain:

\begin{theorem}\label{seg-seg}
Fix an integer $n\geq 2$. If $T: (L^0)^n\to (L^0)^n $ is a bijection which is local and maps each $L^0$-line segment onto an $L^0$-line segment, then $T$ must be an $L^0$-affine linear mapping.
\end{theorem}

 Theorem \ref{seg-seg} will be used in our forthcoming study to give representations of fully order preserving and fully order reversing operators acting on the set of $L^0$-lower semi-continuous $L^0$-convex functions on $(L^0)^n$, and on general complete random normed modules.

\vspace{0.5cm}

{\bf Acknowledgments.}  The first author was supported by the Natural Science Foundation of China (Grant No.11701531) and the Fundamental Research Funds for the Central Universities, China University of Geosciences (Wuhan) (Grant No. CUGL170820). The second author was supported by the Natural Science Foundation of China(Grant No.11501580).


\begin{thebibliography}{1}

\bibitem{AM} S. Artstein-Avidan, V. Milman, The concept of duality
in convex analysis, and the characterization of the Legendre
transform, Ann. Math., 2009, 169(2): 661-674.

\bibitem{AS} S. Artstein-Avidan, B.A.Slomka, The fundamental theorems of affine and projective
geometry revisited, Commun. Contemp. Math., 19(5), 1650059 (2017).


\bibitem{CKV} P. Cheridito, M. Kupper, N. Vogelpoth, Conditional analysis on $R^d$, in: Set Optimization and Applications - State of the Art, Springer, 2015, pp. 179-211.

 \bibitem{Cohn} P. M. Cohn, Some remarks on the invariant basis property, Topology, 5(1966) 215-228.


\bibitem{DKK} S. Drapeau, M. Karliczek, M. Kupper, M. Streckfu{\ss}, Brouwer fixed point theorem in $(L^0)^d$, Fixed Point Theory Appl.,
301(1) (2013).

\bibitem{GZ} T. X. Guo, X. L. Zeng, An $L^0({\mathcal F},R)$-valued function's intermediate value theorem and its applications to
random uniform convexity. Acta Math. Sin. (Engl. Ser.) 28(5)(2012) 909-924.

\bibitem{Kar}
C.Kardaras, Uniform integrability and local convexity in $L^{0}$, J. Funct. Anal., 266 (2014) 1913-1927.

\bibitem{KZ}
C.Kardaras, G. \v{Z}itkovi\'{c}, Forward--convex convergence in probability of sequences of nonnegative random variables, Proc. Amer. Math. Soc., 141(3) (2013) 919-929.

\bibitem{KL} T. G. Kvirikashvili, A. A. Lashkhi, Geometrical maps in ring affine geometries, J. Math. Sci., 186(5) (2012) 759-765.

\bibitem{Wu} M. Z. Wu, Farkas' lemma in random locally convex modules and Minkowski-Weyl
type results in $L^0({\mathcal F}, R^n)$, J. Math. Anal. Appl., 404(2)(2013)  300-309.

\bibitem{Zit}
G.\v{Z}itkovi\'{c}, Convex compactness and its applications, Math. Financ. Econ. 3(1) (2010) 1-12.


\end{thebibliography}
\end{document}